\documentclass[10.5pt,letterpaper]{article}

\usepackage[T1]{fontenc}
\usepackage[utf8]{inputenc}
\usepackage{lmodern}
\usepackage[margin=1in]{geometry}
\usepackage{amsmath,amssymb,amsthm,mathrsfs,bm}
\usepackage{microtype}
\usepackage[numbers,sort&compress]{natbib}
\usepackage{xurl}
\usepackage[hidelinks]{hyperref}

\newtheorem{theorem}{Theorem}
\newtheorem{lemma}{Lemma}
\newtheorem{proposition}{Proposition}
\newtheorem{corollary}{Corollary}
\theoremstyle{definition}
\newtheorem{definition}{Definition}

\newif\ifanonymous
\anonymousfalse

\title{Randomized Parallel Complexity of\\ Weakly Convex Optimization}
\ifanonymous
  \author{Anonymous authors}
  \hypersetup{pdfauthor={}}
\else
  \author{Kang Liu\thanks{School of Future Technology, Xi'an Jiaotong University.
    Email: \href{mailto:kanyo@foxmail.com}{\nolinkurl{kanyo@foxmail.com}}.}
  \hypersetup{pdfauthor={Kang Liu}}
\fi
\date{}
\hypersetup{
  pdftitle={Randomized Parallel Complexity of Weakly Convex Optimization},
  pdfsubject={Weakly convex optimization and parallel oracle complexity}
}
}
\begin{document}
\maketitle

\begin{abstract}
We study the minimum number of query rounds needed to find approximate stationary points of nonsmooth weakly convex functions, measured by the gradient of the Moreau envelope. Each round permits multiple simultaneous queries. In the stated parameter regime, the standard subgradient rate is optimal in sufficiently high dimensions, even for randomized algorithms making polynomially many parallel queries per round. We then show how the order structure of the real line enables faster search: a deterministic algorithm using only function values achieves logarithmic round complexity, with a further reduction through parallel queries. Matching lower bounds hold even when each query returns the entire subdifferential. Together, these results clarify the role of dimension in determining the benefits of parallel queries for weakly convex optimization.
\end{abstract}

\noindent\textbf{Keywords:} Weakly convex optimization, oracle complexity, parallel algorithms, randomized algorithms, Moreau envelope.

\section{Introduction}
\label{sec:introduction}

A function $f:\mathbb{R}^d\to\mathbb{R}$ is $\rho$-weakly convex (WC) if $x\mapsto f(x)+(\rho/2)\|x\|^2$ is convex, where $\rho>0$ and $\|\cdot\|$ denotes the Euclidean norm. Weak convexity allows both nonconvexity and nonsmoothness while controlling negative curvature \citep{Vial1983}. First-order methods seek approximate stationary points through successive queries of function values and subgradients. A basic question is how many query rounds are necessary to attain a prescribed stationarity accuracy, and whether randomization or parallel queries can reduce this number.

\textbf{Problem}~~
Let $\mathcal F_d(G,\rho,\Delta)$ denote the class of $\rho$-WC functions satisfying
\[
|f(x)-f(y)|\le G\|x-y\|
~\text{for all }x,y\in\mathbb R^d,
~~
f(0)-\inf f\le\Delta,
\]
where $G>0$ bounds the global Lipschitz constant and $\Delta>0$ bounds the initial objective gap, with the initial point taken to be the origin. We use the subdifferential
$\partial f(x):=\partial_{\mathrm{cvx}}(f+(\rho/2)\|\cdot\|^2)(x)-\rho x$,
where $\partial_{\mathrm{cvx}}$ denotes the convex subdifferential.

Approximate stationarity is measured by the gradient of the Moreau envelope (ME) \citep{DavisDrusvyatskiy2018,DavisDrusvyatskiy2019}, defined by
$f_{1/(2\rho)}(x):=\min_{y\in\mathbb R^d}\{f(y)+\rho\|y-x\|^2\}$.
For $\varepsilon>0$, the condition $\|\nabla f_{1/(2\rho)}(x)\|\le\varepsilon$ guarantees that some point within distance $\varepsilon/(2\rho)$ of $x$ has a subgradient of norm at most $\varepsilon$.

Each query returns the exact function value and the entire subdifferential, $\mathcal O_f(x):=(f(x),\partial f(x))$. For an integer $P\ge1$, the algorithm submits at most $P$ queries in each round, chosen jointly using its internal randomness and replies from earlier rounds. Replies from the current round become available only after all its queries have been chosen. Thus, $P=1$ gives the sequential model, while larger $P$ permits parallel evaluations \citep{DiakonikolasGuzman2020}. Counting rounds measures the sequential interaction required with the oracle. Query points, output points, and computation between rounds are unrestricted.

We seek the smallest integer $T\ge0$ for which there exists a randomized algorithm using at most $P$ queries per round and at most $T$ rounds on every execution, such that
\[
\Pr_{\mathcal A}\!\left\{
\|\nabla f_{1/(2\rho)}(\widehat x)\|\le\varepsilon
\right\}\ge\frac34
~~
\text{for every }f\in\mathcal F_d(G,\rho,\Delta).
\]
Here, $\widehat x$ is the algorithm's output, and the probability is taken over its internal randomness. A hard instance may depend on the algorithm but is independent of its realized random seed.

\textbf{Results}~~
We construct high-dimensional hard instances to establish a round lower bound for randomized parallel algorithms. Suppose $\varepsilon^2\le c_0\min\{G^2,\rho\Delta\}$ and $\Delta\le G^4/(\rho\varepsilon^2)$. For every fixed $k\ge0$, there exists a sufficiently large dimension $d$ in which the optimal round complexity is $\Theta(\rho G^2\Delta/\varepsilon^4)$ simultaneously for all $1\le P\le d^k$. This lower bound matches the standard subgradient upper bound \citep{DavisDrusvyatskiy2018}, showing that the standard rate remains optimal even with randomization and polynomially many parallel queries per round.

We then exploit the order structure of the real line through a deterministic interval search using only function values. Parallel queries allow finer subdivisions of the search interval in each round. If $\varepsilon\le c_1G$ and $\rho\Delta\ge C_1\varepsilon^2$, this algorithm attains the optimal round complexity $\Theta(1+\log(\rho\Delta/\varepsilon^2)/\log(P+1))$ for every $P\ge1$. The matching lower bound holds for randomized algorithms receiving the entire subdifferential. Thus, the order structure enables logarithmic sequential complexity, with a further reduction through parallel queries. Here, $c_0,c_1>0$ and $C_1>1$ are universal constants.

\textbf{Roadmap}~~
Section~\ref{sec:related-work} reviews related work. Sections~\ref{sec:construction} and~\ref{sec:randomized-complexity} develop the construction and establish the high-dimensional complexity bounds, respectively. Section~\ref{sec:one-dimensional} derives the matching one-dimensional bounds. The main text presents the central construction and complexity arguments; the appendices provide detailed regularity estimates, auxiliary calculations, and supporting proofs.

\section{Related work}
\label{sec:related-work}

The first-order theory of weakly convex optimization has progressed from asymptotic convergence to finite-time stationarity guarantees. Early work established asymptotic convergence for stochastic subgradient and composite methods \citep{DuchiRuan2018}, while proximal regularization yielded nonasymptotic guarantees \citep{DavisGrimmer2019}. Using the gradient of the Moreau envelope to measure stationarity, analyses of the subgradient method and more general model-based methods established the complexity bound $O(\rho G^2\Delta/\varepsilon^4)$ \citep{DavisDrusvyatskiy2018,DavisDrusvyatskiy2019}. These upper bounds raise the question of whether this query complexity is intrinsic to the function class. Under an oracle returning exact function values and the entire subdifferential, \citet{LiPan2026} proved a matching deterministic lower bound for $\Delta\le G^2/\rho$, establishing the optimal deterministic query complexity in the corresponding parameter regime.

This deterministic optimality leaves open whether randomization and parallel queries can improve the complexity. Assessing parallelism requires measuring the number of sequential query rounds rather than the total number of queries. This perspective was developed in parallel oracle models for convex optimization \citep{Nemirovski1994}, with subsequent lower bounds covering randomized algorithms and several smoothness classes \citep{DiakonikolasGuzman2020}. For nonconvex stationarity, hard instances can force algorithms to acquire information about new directions sequentially; this mechanism underlies randomized query lower bounds for smooth nonconvex optimization \citep{CarmonEtAl2020}. Parallel lower bounds further show that, in sufficiently high dimensions, polynomially many queries per round do not improve the optimal complexity order for the corresponding smooth function classes \citep{ZhouEtAl2025}.

Dimension also affects the complexity of stationarity search. For smooth optimization over compact domains, gradient-flow trapping demonstrates the advantage of low-dimensional search over general descent methods \citep{BubeckMikulincer2020}. In one-dimensional unconstrained optimization, function values suffice for optimal logarithmic search complexity under smoothness \citep{ChewiEtAl2023}, with extensions to higher-order H\"older smoothness \citep{ZhouSuzuki2026}. Low-dimensional improvements also occur in nonsmooth optimization. Using Goldstein stationarity, which measures the norm of convex combinations of subgradients at nearby points, \citet{DavisEtAl2022} obtain improved complexity bounds for Lipschitz and weakly convex functions through a cutting-plane method. Under the same stationarity criterion, \citet{KornowskiShamir2022} establish a randomized one-dimensional lower bound logarithmic in the inverse neighborhood radius, using binary strings encoded by nested intervals.

We establish matching randomized parallel bounds for nonsmooth weakly convex optimization under the Moreau-envelope stationarity criterion. In high dimensions, the standard subgradient rate remains optimal with polynomially many queries per round, and the matching range extends to $\Delta\le G^4/(\rho\varepsilon^2)$. We then characterize the faster search enabled by the order structure of the real line, including its dependence on the parallel query budget. Both lower bounds allow access to the entire subdifferential, while the one-dimensional upper bound requires only function values.

\section{Construction of high-dimensional hard instances}
\label{sec:construction}

To prove a lower bound on query rounds, we relate the information supplied by the oracle to approximate stationarity. We construct functions with two properties: exact oracle locality restricts how replies depend on successive coordinate groups, and a lower bound on the ME gradient excludes approximate stationarity when the coordinates in the final block are sufficiently small.

Starting from the block coupling of \citet{LiPan2026}, we build $K$ component functions, called stages, on disjoint coordinate sets. A parameter $N$ controls the size of each stage. For $K\le N$, the final function has Lipschitz and weak convexity parameters at most one, an initial objective gap of order $K$, and order $KN^2$ coordinate groups. When the final-block projections are sufficiently small, its ME gradient norm is bounded below by a constant multiple of $N^{-1/2}$. The construction therefore allows the initial gap and the number of coordinate groups to grow together while maintaining the required stationarity bound.

Let $F:\mathbb{R}^n\to\mathbb{R}$ be a finite-valued WC function, and write $[n]:=\{1,\ldots,n\}$. For a nonempty set $C$, let $\operatorname{dist}(0,C):=\inf_{g\in C}\|g\|$. Consider nested coordinate index sets $\varnothing=\mathcal{I}_{-1}\subseteq\mathcal{I}_0\subseteq\cdots\subseteq\mathcal{I}_H=[n]$, where $\mathcal{I}_t\setminus\mathcal{I}_{t-1}$ is the $t$th group of coordinate indices. Let $\Pi_t$ be the coordinate projection retaining the coordinates in $\mathcal{I}_t$ and setting the others to zero.

\begin{definition}[Exact oracle locality]
\label{def:oracle-locality}
Given $a>0$, the function $F$ has exact oracle locality at scale $a$ with respect to $(\mathcal{I}_t)_{t=-1}^H$ if, for every $t\in\{0,\ldots,H\}$ and every $z$ satisfying $|z_j|\le a$ for $j\notin\mathcal{I}_{t-1}$,
\(
F(z)=F(\Pi_tz),
~
\partial F(z)=\partial F(\Pi_tz)
\subseteq\operatorname{range}\Pi_t.
\)
\end{definition}

Thus, when the coordinates outside $\mathcal{I}_{t-1}$ are sufficiently small, the oracle reply may depend on the new coordinates in $\mathcal{I}_t\setminus\mathcal{I}_{t-1}$, but coordinates outside $\mathcal{I}_t$ affect neither the function value nor the entire subdifferential.

\subsection{The base function and its oracle locality}
\label{subsec:base-function}

Fix integers $N\ge60$ and $2\le M\le N/12$. Write $u=(u^{(1)},\ldots,u^{(M)})\in\mathbb{R}^{MN}$, with $u^{(r)}=(u^{(r)}_1,\ldots,u^{(r)}_N)$. Coordinates are indexed by pairs $(r,j)$, where $r$ identifies the block and $j$ the coordinate within it; the corresponding linear index is $(r-1)N+j$.

Set $a:=1/(8N^2)$, $L:=\lfloor3N/16\rfloor$, and $H_b:=(M-1)L+N-1$. For $0\le t\le H_b$, define
\[
\mathcal{J}_t
:=
\{(r,j)\in[M]\times[N]:(r-1)L+j-1\le t\},
~~
\mathcal{J}_{-1}:=\varnothing.
\]
Let $Q_t$ be the orthogonal projection retaining the coordinates in $\mathcal{J}_t$. The first coordinate of block $r$ enters at $t=(r-1)L$, and $\mathcal{J}_{H_b}=[M]\times[N]$.

The base function combines affine branches within each block with progress terms linking adjacent blocks. Set $c_j:=1+(N-j)/(16N)$ and define $q(s):=0$ for $s\le0$, $q(s):=3s^2-2s^3$ for $0<s<1$, and $q(s):=1$ for $s\ge1$. For $v\in\mathbb{R}^N$, let
\[
\phi(v):=
\max\left\{0,\max_{1\le j\le N}\left(c_j-\frac N8v_j\right)\right\},
~~
\psi(v):=
\frac1N\sum_{j=1}^N q\!\left(\frac N4(v_j-2a)\right).
\]
The intercept gap $c_j-c_{j+1}=1/(16N)$ separates the affine branches when the corresponding coordinates are small. The function $\psi$ takes values in $[0,1]$ and measures progress within a block. The shift by $2a$ ensures that each summand and its derivative vanish whenever $|v_j|\le a$.

Let $p_0:=1$, $p_r:=\psi(u^{(r)})$, $\sigma_r:=\phi(u^{(r)})+p_{r-1}-5/4$, and $\sigma(u):=(\sigma_1,\ldots,\sigma_M)$. Define
\[
S(u):=\|[\sigma(u)]_+\|,
~~
h_{M,N}(u):=\frac2N\left(S(u)-\sum_{r=1}^M p_r\right),
\]
where $[w]_+$ denotes the componentwise positive part of $w$. The negative sum decreases as the progress values grow. Increasing $p_{r-1}$ raises $\sigma_r$, while reducing $\phi(u^{(r)})$ offsets this increase. This couples progress in adjacent blocks.

\begin{proposition}[Properties of the base function]
\label{prop:base-properties}
The function $h_{M,N}$ is $3/4$-WC, with $\operatorname{Lip}(h_{M,N})\le1/4+(3/4)\sqrt{M/N}<1/2$. It satisfies $h_{M,N}\ge-2M/N$, $h_{M,N}(0)\in[0,2/N)$, and $h_{M,N}(0)-\inf h_{M,N}\le3M/N$. All its subgradients are componentwise nonpositive. If $S(u)>0$, then $\|g\|\ge1/(4\sqrt N)$ for every $g\in\partial h_{M,N}(u)$. Moreover, $p_M<27/32$ implies $S(u)>0$.
\end{proposition}

In particular, insufficient progress in the last block guarantees the subgradient norm lower bound. The proof is given in Appendix~\ref{app:base-properties}.

\begin{lemma}[Exact oracle locality of the base function]
\label{lem:base-locality}
For every integer $0\le t\le H_b$, if $|u^{(r)}_j|\le a$ for all $(r,j)\notin\mathcal{J}_{t-1}$, then
\(
h_{M,N}(u)=h_{M,N}(Q_tu),
~
\partial h_{M,N}(u)=\partial h_{M,N}(Q_tu)
\subseteq\operatorname{range}Q_t.
\)
\end{lemma}

The proof establishes $h_{M,N}=h_{M,N}\circ Q_t$ on open neighborhoods of both $u$ and $Q_tu$. On these neighborhoods, omitted coordinates affect neither the progress terms nor the active block terms. The local identity yields equality of the full subdifferentials. The details are given in Appendix~\ref{app:base-locality}.

\subsection{Construction of a single stage}
\label{subsec:single-stage}

We will combine shifted copies of a single-stage function by taking their maximum, which preserves common Lipschitz and weak convexity bounds. Subgradients of the maximum are convex combinations of subgradients of the active branches, namely those attaining the maximum. To control the resulting decrease in norm, we modify the base function so that subgradients of norm below $1/4$ occur only within a fixed interval of function values.

For the stage construction, take $N\ge120$, $m:=\lfloor N/12\rfloor$, $d_0:=mN$, and $h_N:=h_{m,N}$. Let $e:=(0,\ldots,0,\mathbf{1}_N/\sqrt N)\in\mathbb{R}^{d_0}$, where $\mathbf{1}_N$ is the all-ones vector; thus, $e$ is a unit vector supported on block $m$. Define $\tau_N(s):=\int_0^{(s-1/\sqrt N)_+}\min\{t/20,1/4\}\,dt$ and
\[
b_N(u)
:=
h_N(u)+\frac14\|[-u-2a\mathbf{1}]_+\|
-\tau_N(\langle e,u\rangle).
\]
The penalty on negative coordinates controls subgradient norms at large function values. The terminal descent term becomes active when $\langle e,u\rangle>1/\sqrt N$, allowing further decrease and controlling subgradient norms at small function values.

\begin{lemma}[Single-stage properties]
\label{lem:single-stage}
The function $b_N$ is $4/5$-WC and globally Lipschitz with constant less than $31/32$. For every $u$ and $g\in\partial b_N(u)$, the vector $g$ is componentwise nonpositive, $\|g\|\ge1/(4\sqrt N)$, and $\langle g,u\rangle\ge b_N(u)-2$. If $b_N(u)\notin(-1,1/4)$, then $\|g\|\ge1/4$. Moreover, $b_N(0)=h_N(0)\in[0,2/N)$, and $\langle e,u\rangle\le1/\sqrt N$ implies $b_N(u)\ge-1/6$. The function $b_N$ is unbounded below along the positive $e$ direction.
\end{lemma}

The three terms defining $b_N$ have componentwise nonpositive subgradients, so their contributions cannot cancel. At function values at least $1/4$, the penalty has gradient norm $1/4$; at values at most $-1$, the terminal descent term has slope magnitude $1/4$. The full proof, including the regularity and inner-product estimates, is given in Appendix~\ref{app:single-stage}.

\subsection{Combining stages while preserving the subgradient bound}
\label{subsec:stage-combination}

We choose shifts that separate the intervals of function values where small subgradients can occur. At most one active stage can then have a small subgradient, while the remaining active stages satisfy a stronger bound. The following lemma quantifies the resulting bound for their maximum.

\begin{lemma}[Combination with separated value intervals]
\label{lem:separated-intervals}
Let $b_i:\mathbb{R}^{d_i}\to\mathbb{R}$, $i=1,\ldots,K$, be finite-valued WC functions on disjoint coordinate blocks. Suppose there are a bounded interval $I$ and constants $0<\beta\le\gamma$ such that $\operatorname{dist}(0,\partial b_i(v))\ge\beta$ for every $i,v$, and $\operatorname{dist}(0,\partial b_i(v))\ge\gamma$ whenever $b_i(v)\notin I$. Choose shifts satisfying $|s_i-s_j|>\operatorname{length}(I)$ for $i\ne j$, and let
\[
\mathcal{M}(z)
:=
\max\left\{
\max_{1\le i\le K}\bigl(b_i(z^{(i)})-s_i\bigr),\,\nu
\right\},
\]
where $\nu$ is constant. If $\mathcal{M}(z)>\nu$, then every $g\in\partial\mathcal{M}(z)$ satisfies $\|g\|\ge(\beta^{-2}+(K-1)\gamma^{-2})^{-1/2}$. If the stages share a Lipschitz bound and a weak convexity parameter, then $\mathcal{M}$ satisfies the same bounds.
\end{lemma}

\begin{proof}
Let $\mathcal{A}:=\{i:b_i(z^{(i)})-s_i=\mathcal{M}(z)\}$ be the set of active stages. Since the constant branch is inactive, every $g\in\partial\mathcal{M}(z)$ has a representation $g=\sum_{i\in\mathcal{A}}\lambda_i\widetilde g_i$, where $\lambda_i\ge0$, $\sum_{i\in\mathcal{A}}\lambda_i=1$, $g_i\in\partial b_i(z^{(i)})$, and $\widetilde g_i$ is the embedding of $g_i$ into its coordinate block.

Each active stage satisfies $b_i(z^{(i)})=\mathcal{M}(z)+s_i$. The separation of the shifts implies that at most one active stage has its function value in $I$. Since distinct stages have orthogonal coordinate supports, Cauchy--Schwarz gives
\[
\|g\|^2
=
\sum_{i\in\mathcal{A}}\lambda_i^2\|g_i\|^2
\ge
\frac1{\sum_{i\in\mathcal{A}}\|g_i\|^{-2}}
\ge
\frac1{\beta^{-2}+(K-1)\gamma^{-2}}.
\]
The regularity claims follow from the Lipschitz property of the maximum and a common convexifying quadratic.
\end{proof}

Fix an integer $1\le K\le N$. Write $z=(z^{(1)},\ldots,z^{(K)})\in\mathbb{R}^n$, where $z^{(i)}\in\mathbb{R}^{d_0}$ and $n:=Kd_0$, and define
\[
F_K(z)
:=
\max\left\{
\max_{1\le i\le K}\bigl(b_N(z^{(i)})-2(i-1)\bigr),\,-2K
\right\}.
\]
Adjacent shifts differ by $2$, exceeding the length of $(-1,1/4)$. The constant branch bounds $F_K$ below by $-2K$. Each stage can decrease independently along its terminal direction, so $\inf F_K=-2K$ and $F_K(0)=h_N(0)$, giving an initial objective gap of order $K$.

Applying Lemma~\ref{lem:separated-intervals} with $\beta=1/(4\sqrt N)$ and $\gamma=1/4$ gives
\[
F_K(z)>-2K
~\Longrightarrow~
\operatorname{dist}(0,\partial F_K(z))\ge\sigma_K,
~~
\sigma_K:=\frac1{4\sqrt{N+K-1}}.
\]
Since $K\le N$, we have $\sigma_K\ge1/(4\sqrt{2N})$, so the combined function retains a subgradient lower bound of order $N^{-1/2}$.

The function $F_K$ is $4/5$-WC and has Lipschitz constant less than $31/32$. We will also use the inner-product estimate $\langle g,z\rangle\ge F_K(z)-2\ge-2K-2$ for every $g\in\partial F_K(z)$. Indeed, if stage $i$ is active, then $b_N(z^{(i)})=F_K(z)+2(i-1)$, and Lemma~\ref{lem:single-stage} bounds the corresponding inner product below by $F_K(z)-2$. When the constant branch is active, its zero subgradient also satisfies this bound. Taking convex combinations proves the bound.

\subsection{Oracle locality of the combined function}
\label{subsec:combined-locality}

We now establish exact oracle locality for $F_K$ with coordinate groups ordered stage by stage. Index coordinates by triples $(i,r,j)$, where $i$ identifies the stage, $r$ the block within that stage, and $j$ the coordinate within the block. The corresponding linear index is $(i-1)mN+(r-1)N+j$.

Set $J:=(m-1)L$ and $H:=K(J+1)-1$. The terminal descent term depends on the entire last block, so its coordinates enter together in the final group of each stage. We retain the base ordering for indices below $J$ and place all remaining coordinates in that final group. For $0\le t\le H$, define
\[
\mathcal{I}_t
:=
\bigl\{(i,r,j)\in[K]\times[m]\times[N]:
(i-1)(J+1)
+\min\{(r-1)L+j-1,J\}\le t\bigr\},
\]
and let $\mathcal{I}_{-1}:=\varnothing$. Let $\Pi_t$ be the orthogonal projection retaining the coordinates in $\mathcal{I}_t$.

By construction, $\mathcal{I}_H=[K]\times[m]\times[N]$ and $\{(K,m,j):j\in[N]\}\subseteq\mathcal{I}_H\setminus\mathcal{I}_{H-1}$. Thus, the entire last block of the final stage belongs to the last coordinate group. There are $H+1$ groups, with $H\ge KJ\ge KN^2/288$.

\begin{lemma}[Exact oracle locality of the combined function]
\label{lem:combined-locality}
For every integer $0\le t\le H$, if $|z_{i,r,j}|\le a$ for all $(i,r,j)\notin\mathcal{I}_{t-1}$, then
\(
F_K(z)=F_K(\Pi_tz),
~
\partial F_K(z)=\partial F_K(\Pi_tz)
\subseteq\operatorname{range}\Pi_t.
\)
\end{lemma}

\begin{proof}
Write $t=(i-1)(J+1)+s$ uniquely with $0\le s\le J$. All coordinates of stages preceding $i$ are retained by $\Pi_t$. The last block of stage $i$ lies outside $\mathcal{I}_{t-1}$, so $\langle e,z^{(i)}\rangle\le\sqrt N\,a<1/\sqrt N$. Its terminal descent term therefore vanishes on a neighborhood, and deleted coordinates lie in the zero region of the penalty. If $s<J$, Lemma~\ref{lem:base-locality} gives the required local identity for stage $i$; if $s=J$, this stage is retained in full.

Lemma~\ref{lem:single-stage} gives $b_N(z^{(i)})\ge-1/6$. For every $k>i$, all coordinates of stage $k$ have magnitude at most $a$, so its progress terms, penalty, and terminal descent term vanish, and $|b_N(z^{(k)})-h_N(0)|\le a/4$. The shifts therefore place every subsequent stage strictly below stage $i$. The constant branch is also strictly below it, since $-1/6-2(i-1)>-2K$.

The preceding stages are retained in full, while the local identity for stage $i$ preserves its value and full subdifferential after projection. The later stages and the constant branch are inactive near both $z$ and $\Pi_tz$. Thus, the active branches have identical values and subdifferentials at the two points, with all subgradients supported on $\mathcal{I}_t$. The subdifferential formula for a maximum proves the claim.
\end{proof}

\subsection{High-dimensional embedding and ME gradient bounds}
\label{subsec:radial-extension}

We embed $F_K$ into a higher-dimensional space and add a radial branch so that oracle replies outside a bounded region are independent of the embedding. The inner-product estimate above controls the subgradients where both branches are active, allowing us to preserve the required norm bound.

Let $D\ge n$ and $U\in\mathbb{R}^{D\times n}$ satisfy $U^\top U=I_n$. Set $A:=2K+4(2K+2)/\sigma_K$ and $B:=65(A+1)$, and define $\Phi_U(x):=\max\{F_K(U^\top x),\,\|x\|-A\}$.

\begin{proposition}[Properties of the radial extension]
\label{prop:radial-extension}
For every such $U$, the function $\Phi_U$ is globally $1$-Lipschitz and $4/5$-WC, with $\Phi_U(0)-\inf\Phi_U=2K+h_N(0)\in[2K,3K]$. If $\|x\|\ge64(A+1)$, then $\mathcal{O}_{\Phi_U}(x)=(\|x\|-A,\{x/\|x\|\})$. This oracle identity also holds after replacing $F_K$ by $F_K\circ\Pi_t$ for any $0\le t\le H$. Moreover, $F_K(U^\top x)>-2K$ implies $\operatorname{dist}(0,\partial\Phi_U(x))>\sigma_K/2$.
\end{proposition}

\begin{proof}
The regularity bounds, objective gap, and oracle identities are verified in Appendix~\ref{app:radial-properties}. Here we prove the subgradient lower bound, which controls the interaction between the two branches.

Assume $F_K(U^\top x)>-2K$. If only one branch is active, the claim follows from the subgradient lower bound for $F_K$ or the unit norm of the radial gradient. Suppose both branches are active. Let $v:=Ug$ with $g\in\partial F_K(U^\top x)$, and set $e_x:=x/\|x\|$. Since the radial branch is active, $\|x\|\ge A-2K=4(2K+2)/\sigma_K$. Using $\langle g,U^\top x\rangle\ge-2K-2$, we obtain
\[
\langle v,e_x\rangle
\ge-\frac{2K+2}{\|x\|}
\ge-\frac{\sigma_K}{4}
\ge-\frac{\|v\|}{4}.
\]
Thus the radial gradient cannot nearly cancel an embedded subgradient. More precisely, for every $0\le\lambda\le1$,
\[
\|\lambda v+(1-\lambda)e_x\|
\ge\sqrt{\frac38}\min\{\|v\|,1\}
\ge\sqrt{\frac38}\,\sigma_K
>\frac{\sigma_K}{2}.
\]
The full subdifferential consists of these convex combinations, proving the bound.
\end{proof}

To bound the ME gradient, we apply the subgradient estimate at the proximal point. A small ME gradient places this point near $x$. If the projections of $x$ onto the last block are also small, the embedded function remains above $-2K$ at the proximal point. Let $u_{K,m,j}$ denote the column of $U$ corresponding to coordinate $(K,m,j)$.

\begin{lemma}[Last-block projections and the ME gradient]
\label{lem:terminal-moreau}
If $|\langle u_{K,m,j},x\rangle|\le a$ for every $j\in[N]$, then $\|\nabla(\Phi_U)_{1/2}(x)\|>\sigma_K/4=1/(16\sqrt{N+K-1})$. The same bound holds whenever $\|x\|\ge B$.
\end{lemma}

\begin{proof}
Suppose $\|\nabla(\Phi_U)_{1/2}(x)\|\le\sigma_K/4$. Let $y:=\arg\min_w\{\Phi_U(w)+\|w-x\|^2\}$ and $v:=\nabla(\Phi_U)_{1/2}(x)$. Proximal optimality gives $v=2(x-y)\in\partial\Phi_U(y)$ and $\|x-y\|\le\sigma_K/8$.

Under the projection condition in the lemma, $\langle e,(U^\top y)^{(K)}\rangle\le\sqrt N\,a+\sigma_K/8<1/\sqrt N$. Lemma~\ref{lem:single-stage} bounds the value of the last stage below by $-1/6$, so $F_K(U^\top y)\ge-1/6-2(K-1)>-2K$. Proposition~\ref{prop:radial-extension} gives $\|v\|>\sigma_K/2$, a contradiction.

If instead $\|x\|\ge B$, the same displacement bound gives $\|y\|\ge B-\sigma_K/8>64(A+1)$. By the radial oracle identity, the subdifferential at $y$ contains only the unit radial vector, so $\|v\|=1$, again a contradiction.
\end{proof}

For $1\le K\le N$, the resulting function $\Phi_U$ has uniformly bounded Lipschitz and weak convexity parameters and an initial objective gap of order $K$. Before embedding, $F_K$ has order $KN^2$ ordered coordinate groups and satisfies exact oracle locality; sufficiently small projections onto the final block imply an ME gradient lower bound of order $N^{-1/2}$ for $\Phi_U$. The radial extension makes oracle replies outside a bounded region independent of the embedding, combining oracle locality, a controlled initial gap, and the required stationarity bound in the same construction.

\section{Randomized parallel complexity in high dimensions}
\label{sec:randomized-complexity}

This section establishes round lower bounds for randomized parallel algorithms in high dimensions. Since queries depend on previous oracle replies, they also depend indirectly on the embedding directions. We control this dependence through an auxiliary execution and use concentration under a random embedding to verify the projection conditions needed to exclude approximate stationarity.
We state the results for a target success probability of $3/4$. The same round complexity orders hold for any fixed success probability $p\in(1/2,1)$, with implicit constants that may depend on $p$.

\subsection{Main theorem and the matching regime}
\label{subsec:high-dimensional-result}

\begin{theorem}[High-dimensional round lower bound]
\label{thm:high-dimensional}
There exist universal constants $c_0,c,C>0$ such that, for all $G,\rho,\Delta,\varepsilon>0$ satisfying $\varepsilon^2\le c_0\min\{G^2,\rho\Delta\}$ and every integer $P\ge1$, there is a dimension $D\le C(G/\varepsilon)^{14}\log(PG/\varepsilon)$ with the following property. For every randomized algorithm using at most $P$ queries per round and at most
\begin{equation}
\label{eq:high-dimensional-budget}
\left\lfloor
c\min\left\{
\frac{\rho G^2\Delta}{\varepsilon^4},
\frac{G^6}{\varepsilon^6}
\right\}
\right\rfloor
\end{equation}
rounds on every execution, there is a fixed instance $f\in\mathcal F_D(G,\rho,\Delta)$ such that its output $\widehat x$ satisfies
\[
\Pr_{\mathcal A}\left\{
\|\nabla f_{1/(2\rho)}(\widehat x)\|>\varepsilon
\right\}\ge\frac34.
\]
The instance may depend on the algorithm but is independent of its realized random seed.

For every fixed $k\ge0$, a dimension $D=O_k((G/\varepsilon)^{14}\log(G/\varepsilon))$ can be chosen so that the same lower bound holds simultaneously for all $1\le P\le D^k$.
\end{theorem}

Combining this result with the standard subgradient upper bound gives the following optimal complexity characterization.

\begin{corollary}[Optimal high-dimensional round complexity]
\label{cor:high-dimensional-optimal}
Under the accuracy condition of Theorem~\ref{thm:high-dimensional}, suppose $\Delta\le G^4/(\rho\varepsilon^2)$. For every fixed integer $P\ge1$, the randomized parallel round complexity, with the worst case taken over dimensions, is $\Theta(\rho G^2\Delta/\varepsilon^4)$. For every fixed $k\ge0$, this order also holds in a single sufficiently large dimension $D=O_k((G/\varepsilon)^{14}\log(G/\varepsilon))$, simultaneously for all $1\le P\le D^k$.
\end{corollary}

The condition $\Delta\le G^4/(\rho\varepsilon^2)$ makes the first term in the minimum in Theorem~\ref{thm:high-dimensional} no larger than the second, yielding a lower bound that matches the subgradient upper bound. For fixed $G,\rho,\Delta$, this matching condition holds whenever the target accuracy $\varepsilon$ is sufficiently small. Beyond this range, the theorem still provides a lower bound, but its order remains $G^6/\varepsilon^6$ and no longer grows with $\Delta$.

\subsection{Information revealed by randomized parallel queries}
\label{subsec:randomized-coupling}

We use the function $\Phi_U$, index sets $\mathcal I_t$, projections $\Pi_t$, and parameters $n,a,H,A,B,\sigma_K$ from Section~\ref{sec:construction}. Identify the triple coordinate indices with $j\in[n]$, and let $u_j$ denote the corresponding column of $U$. Draw $U$ uniformly from the matrices in $\mathbb R^{D\times n}$ satisfying $U^\top U=I_n$, independently of the algorithm's internal randomness.

\begin{proposition}[Normalized lower bound for the stage construction]
\label{prop:normalized-stage-bound}
Let $N\ge120$, $1\le K\le N$, and $P\ge1$ be integers. If
\[
D\ge n+2+\frac{2B^2}{a^2}\log\!\bigl(8n(PH+1)\bigr),
\]
then every randomized algorithm using at most $P$ queries per round and at most $H$ rounds on every execution satisfies
\[
\Pr_{U,\mathcal A}\left\{
\|\nabla(\Phi_U)_{1/2}(\widehat x)\|>\frac{\sigma_K}{4}
\right\}\ge\frac34.
\]
Consequently, for every such algorithm, there is a fixed embedding $U_0$ for which the failure probability on $\Phi_{U_0}$ is at least $3/4$.
\end{proposition}

\begin{proof}
Fix an algorithm with a round cap $T\le H$, and first fix its internal random seed. Its queries and output are then measurable functions of the replies received.

To control the dependence of the queries on the random embedding, we introduce an auxiliary execution whose replies in round $t$ depend only on columns indexed by $\mathcal I_t$. This allows us to apply conditional spherical concentration to columns that have not yet been used when the queries are selected.

Define the filtration $\mathscr F_t:=\sigma\{u_j:j\in\mathcal I_{t-1}\}$ for $t=0,\ldots,T$. It records the columns used to determine the auxiliary replies before round $t$. In the auxiliary execution, round $t$, indexed from $0$ to $T-1$, returns the exact function value and entire subdifferential of $\Phi_{U,t}(x):=\max\{F_K(\Pi_tU^\top x),\,\|x\|-A\}$. This function and its replies depend only on columns indexed by $\mathcal I_t$.

All queries in a round are selected before any reply from that round is received. Inductively, every query in round $t$ is therefore $\mathscr F_t$-measurable, the replies are $\mathscr F_{t+1}$-measurable, and the final output $x_{\mathrm{out}}$ is $\mathscr F_T$-measurable. An algorithm that stops early can be padded to $T$ rounds without changing its output.

Conditioned on $\mathscr F_t$, each column $u_j$ with $j\notin\mathcal I_{t-1}$ is uniformly distributed on the unit sphere in the orthogonal complement of the columns indexed by $\mathcal I_{t-1}$. This space has dimension $D-|\mathcal I_{t-1}|\ge D-n$. Let $x_{t,p}$ denote the $p$th auxiliary query in round $t$. Since $x_{t,p}$ is $\mathscr F_t$-measurable, spherical concentration gives, for every $j\notin\mathcal I_{t-1}$,
\[
\Pr\left\{
\|x_{t,p}\|\le B,\,
|\langle u_j,x_{t,p}\rangle|>a
\,\middle|\,\mathscr F_t
\right\}
\le
2\exp\left(-\frac{(D-n-2)a^2}{2B^2}\right).
\]
The same estimate applies to the output, conditioned on $\mathscr F_T$, for every $j\notin\mathcal I_{T-1}$.

Let $\mathcal G$ be the event that every auxiliary query with norm at most $B$ satisfies $|\langle u_j,x_{t,p}\rangle|\le a$ for all $j\notin\mathcal I_{t-1}$, and that, whenever $\|x_{\mathrm{out}}\|\le B$, the output satisfies $|\langle u_j,x_{\mathrm{out}}\rangle|\le a$ for all $j\notin\mathcal I_{T-1}$. Taking expectations in the conditional bound and applying a union bound over at most $PT$ queries, the final output, and $n$ directions gives
\[
\Pr_U(\mathcal G^c)
\le
2n(PT+1)\exp\left(-\frac{(D-n-2)a^2}{2B^2}\right)
\le\frac14.
\]

On $\mathcal G$, the real and auxiliary executions have identical replies. Indeed, if their previous replies agree, they submit the same queries in the next round. At a query of norm at most $B$, the projection conditions and Lemma~\ref{lem:combined-locality}, together with the subdifferential rule for linear composition, give $\mathcal O_{\Phi_U}(x_{t,p})=\mathcal O_{\Phi_{U,t}}(x_{t,p})$. At a query of norm greater than $B$, both replies are determined by the same radial branch, by Proposition~\ref{prop:radial-extension}, including its assertion for the truncated functions. Induction therefore gives identical transcripts and outputs.

Since $T\le H$, we have $\mathcal I_{T-1}\subseteq\mathcal I_{H-1}$. The entire last block of the final stage lies in $\mathcal I_H\setminus\mathcal I_{H-1}$, so none of its columns is among those generating $\mathscr F_T$. If the common output has norm at most $B$, the event $\mathcal G$ ensures that all its projections onto this block have magnitude at most $a$. If its norm exceeds $B$, the distant-output condition applies instead. In either case, Lemma~\ref{lem:terminal-moreau} gives $\|\nabla(\Phi_U)_{1/2}(x_{\mathrm{out}})\|>\sigma_K/4$.

Thus, for every fixed algorithmic seed, the failure probability over $U$ is at least $3/4$. Averaging over the independent seed proves the joint probability bound. Averaging over $U$ then yields a fixed embedding $U_0$ with algorithmic failure probability at least $3/4$.
\end{proof}

The construction satisfies $n,H=O(KN^2)$ and $B^2/a^2=O(K^2N^5)$. Hence the dimension condition in Proposition~\ref{prop:normalized-stage-bound} can be met with $D=O(K^2N^5\log(PKN))$. Increasing $P$ increases the number of queries covered by the concentration estimate, but does not change which coordinate groups determine an auxiliary reply. It therefore enters this sufficient dimension bound only through a logarithmic factor.

\subsection{Small initial objective gaps}
\label{subsec:small-gap}

The stage construction has a normalized initial objective gap of at least $2$. To cover smaller gaps, we use the base function $h_{M,N}$ directly and vary the number of blocks $M$. Its initial gap is at most $3M/N$, while the round lower bound below is at least $MN/12$ and the ME gradient threshold remains of order $N^{-1/2}$. Thus, once $N$ is chosen according to the target accuracy, varying $M$ controls the gap while retaining a round lower bound that depends linearly on it.

Take integers $N\ge60$ and $2\le M\le N/12$. Set $n_s:=MN$, $s:=\sqrt{M/N}$, $A_s:=34s+2s^2$, and $B_s:=73s$. For $U\in\mathbb R^{D\times n_s}$ satisfying $U^\top U=I_{n_s}$, define $\Psi_U(x):=\max\{h_{M,N}(U^\top x),\,\|x\|-A_s\}$. Proposition~\ref{prop:base-properties} implies that $\Psi_U$ is globally $1$-Lipschitz and $3/4$-WC. The radial branch is inactive at the origin, and $\inf\Psi_U\ge\inf h_{M,N}$, so $\Psi_U(0)-\inf\Psi_U\le3M/N$.

Use the index sets $\mathcal J_t$ and projections $Q_t$ from Section~\ref{subsec:base-function}. For this construction, small projections onto the second half of the last block suffice for the required ME gradient bound. Accordingly, set $H_s:=(M-1)L+\lfloor N/2\rfloor$. Then $H_s\ge MN/12$, and the coordinates in the last block with $j>\lfloor N/2\rfloor$ all lie outside $\mathcal J_{H_s-1}$.

\begin{proposition}[Normalized lower bound for small gaps]
\label{prop:normalized-small-gap}
For the above $M,N$ and any integer $P\ge1$, suppose
\[
D\ge n_s+2+\frac{2B_s^2}{a^2}
\log\!\bigl(8n_s(PH_s+1)\bigr),
~~
a=\frac1{8N^2}.
\]
If $U$ is drawn uniformly from the orthonormal embeddings, independently of the algorithm's internal randomness, then every randomized algorithm using at most $P$ queries per round and at most $H_s$ rounds on every execution satisfies
\[
\Pr_{U,\mathcal A}\left\{
\|\nabla(\Psi_U)_{1/2}(\widehat x)\|>\frac1{16\sqrt N}
\right\}\ge\frac34.
\]
\end{proposition}

\begin{proof}
We use two properties of $\Psi_U$, proved in Appendix~\ref{app:small-gap-construction}. First, for $\|x\|\ge72s$, the radial branch strictly dominates the base-function branch in both $\Psi_U$ and its truncated versions $\max\{h_{M,N}(Q_tU^\top x),\,\|x\|-A_s\}$. Second, writing $u_{M,j}$ for the column of $U$ corresponding to coordinate $j$ of the last block, we have $\|\nabla(\Psi_U)_{1/2}(x)\|>1/(16\sqrt N)$ whenever $|\langle u_{M,j},x\rangle|\le a$ for all $j>\lfloor N/2\rfloor$. The same ME gradient bound holds whenever $\|x\|\ge B_s$.

Apply the auxiliary execution and coupling argument of Proposition~\ref{prop:normalized-stage-bound}, replacing $(\mathcal I_t,\allowbreak\Pi_t,\allowbreak\Phi_U,\allowbreak n,\allowbreak H,\allowbreak B)$ with $(\mathcal J_t,Q_t,\Psi_U,n_s,H_s,B_s)$. Lemma~\ref{lem:base-locality} gives equality of bounded-query replies, and radial dominance gives equality for distant queries.

After at most $H_s$ rounds, the last-block coordinates with $j>\lfloor N/2\rfloor$ remain outside the index set determining the auxiliary output. Under the stated dimension condition, concentration and the union bound ensure that the ME gradient bound applies to the output with probability at least $3/4$.
\end{proof}

Since $B_s^2/a^2=O(MN^3)$, the dimension condition in Proposition~\ref{prop:normalized-small-gap} can be met with $D=O(MN^3\log(PN))$.

\subsection{Scaling and matching upper bounds}
\label{subsec:scaling}

Define the dimensionless parameters $E:=G^2/\varepsilon^2$ and $r:=\rho\Delta/G^2$. Then $\rho G^2\Delta/\varepsilon^4=rE^2$ and $G^6/\varepsilon^6=E^3$.

For a normalized function $\widetilde f$, define $f(x):=(G^2/\rho)\widetilde f((\rho/G)x)$. If $\widetilde f$ has Lipschitz constant and weak convexity parameter at most $1$, then the corresponding parameters of $f$ are at most $G$ and $\rho$. Moreover,
\[
\partial f(x)
=
G\,\partial\widetilde f\!\left(\frac{\rho}{G}x\right),
~~
\nabla f_{1/(2\rho)}(x)
=
G\,\nabla\widetilde f_{1/2}\!\left(\frac{\rho}{G}x\right).
\]
Query locations, function values, and entire subdifferentials can all be transformed by these known scalings. Thus, the transformation preserves both the number of rounds and the number of parallel queries per round. The normalized problem therefore has an initial gap limit of $r$ and a target ME gradient norm of $E^{-1/2}$.

Both constructions have ME gradient thresholds of order $N^{-1/2}$, so we choose $N\asymp E$. For larger normalized gaps, choosing $K\asymp\min\{r,E\}$ gives a round lower bound of order $KN^2\asymp\min\{rE^2,E^3\}$. For smaller gaps, choosing $M\asymp\min\{r,1\}N$ gives the required order $rE^2$. The restriction $K\le N$ preserves the ME gradient threshold of order $N^{-1/2}$ in the stage construction and causes its round lower bound to saturate at order $E^3$.

\begin{proof}[Proof of Theorem~\ref{thm:high-dimensional}]
Choose $c_0\le1/122880$. The accuracy condition implies $E\ge122880$ and $rE\ge122880$.

First suppose $r\ge8$. Set $N:=\lfloor E/1024\rfloor$ and $K:=\min\{\lfloor r/4\rfloor,N\}$, and use $\widetilde f:=\Phi_U$. Then $N\ge120$, $1\le K\le N$, and $f(0)-\inf f\le3KG^2/\rho\le\Delta$. The scaled ME gradient threshold satisfies
\[
\frac{G}{16\sqrt{N+K-1}}
\ge
\frac{G}{16\sqrt{2N}}
\ge
\sqrt2\,\varepsilon.
\]
Since $N=\Theta(E)$ and $K=\Theta(\min\{r,E\})$, Proposition~\ref{prop:normalized-stage-bound} gives
\[
H\ge\frac{KN^2}{288}
=
\Omega\!\left(\min\{rE^2,E^3\}\right).
\]

Next suppose $r<8$. Set $r':=\min\{r,1\}$, $N:=\lfloor E/256\rfloor$, and $M:=\lfloor r'N/12\rfloor$, and use $\widetilde f:=\Psi_U$. The accuracy condition ensures $N\ge60$, $2\le M\le N/12$, $N\ge E/512$, and $M\ge r'N/24$. The objective gap satisfies
\[
f(0)-\inf f
\le
\frac{G^2}{\rho}\frac{3M}{N}
\le
\frac{r'G^2}{4\rho}
\le\Delta,
\]
and the scaled ME gradient threshold is $G/(16\sqrt N)\ge\varepsilon$. Proposition~\ref{prop:normalized-small-gap} therefore gives
\[
H_s
\ge
\frac{MN}{12}
\ge
\frac{r'N^2}{288}
=
\Omega(r'E^2)
=
\Omega(rE^2),
\]
where the last equality uses $r'\ge r/8$. Since $r<E$ in this case, this is the required order $\Omega(\min\{rE^2,E^3\})$.

The sufficient dimensions for the two constructions satisfy, respectively,
\[
\begin{aligned}
D
=
O\!\left(K^2N^5\log(PKN)\right)
=
O\!\left(E^7\log(P\sqrt E)\right),~~
D
=
O\!\left(MN^3\log(PN)\right)
=
O\!\left(r'E^4\log(P\sqrt E)\right).
\end{aligned}
\]
Both imply the stated dimension bound.

For polynomial batch widths, choose $D:=\lceil C_kE^7\log E\rceil$. The intrinsic dimensions and round caps of both constructions are $O(E^3)$. For $P\le D^k$, we have $\log(PE)=O_k(\log E)$, so a sufficiently large $C_k$ satisfies the concentration conditions simultaneously for all such $P$.

Finally, choose the universal constant $c>0$ small enough that the budget in \eqref{eq:high-dimensional-budget} does not exceed the applicable cap $H$ or $H_s$. The two normalized propositions give joint failure probability at least $3/4$. Averaging over the embedding yields a fixed instance with the same algorithmic failure bound.
\end{proof}

When $r\le1$, equivalently $\Delta\le G^2/\rho$, the small-gap construction also gives the sharper sufficient dimension $D=O((\rho\Delta G^6/\varepsilon^8)\log(PG/\varepsilon))$.

\begin{proof}[Proof of Corollary~\ref{cor:high-dimensional-optimal}]
For any $f\in\mathcal F_D(G,\rho,\Delta)$, start from $x_0:=0$, choose $\eta:=\sqrt{\Delta/(\rho G^2T)}$, and perform $T$ subgradient updates $x_{t+1}:=x_t-\eta g_t$ with $g_t\in\partial f(x_t)$. Return $x_R$, where $R$ is uniform on $\{0,\ldots,T-1\}$ and independent of the trajectory. The standard ME descent analysis \citep{DavisDrusvyatskiy2018,DavisDrusvyatskiy2019} gives
\[
\mathbb E\|\nabla f_{1/(2\rho)}(x_R)\|^2
\le
8\sqrt{\frac{\rho G^2\Delta}{T}}.
\]
For $T\ge1024\rho G^2\Delta/\varepsilon^4$, the right-hand side is at most $\varepsilon^2/4$. Markov's inequality gives success probability at least $3/4$, using one query per round. More generally, taking $T\ge64(1-p)^{-2}\rho G^2\Delta/\varepsilon^4$ guarantees success probability at least $p$ for any fixed $p\in(1/2,1)$.

The condition $\Delta\le G^4/(\rho\varepsilon^2)$ is equivalent to $r\le E$. In this regime, Theorem~\ref{thm:high-dimensional} gives the matching lower bound $\Omega(rE^2)=\Omega(\rho G^2\Delta/\varepsilon^4)$. This lower bound also applies to every target success probability $p>1/2$, since the theorem gives a fixed instance on which the algorithm succeeds with probability at most $1/4$. Its simultaneous dimension guarantee yields the conclusion for polynomial batch widths.
\end{proof}

The coupling argument converts the oracle locality and stationarity properties of the construction into a round lower bound for randomized parallel algorithms. After scaling, this bound matches the standard subgradient upper bound when $\Delta\le G^4/(\rho\varepsilon^2)$, under the stated accuracy condition. Thus, in this regime, polynomially many parallel queries per round do not improve the worst-case rate in sufficiently high dimensions, even for randomized algorithms.

\section{Optimal round complexity in one dimension}
\label{sec:one-dimensional}

We next examine how the order structure of the real line permits faster stationarity search. Function-value comparisons allow an algorithm to shrink a search interval while ensuring that a regularized objective has an interior minimizer on the current interval. Parallel queries allow finer subdivisions in each round, and weak convexity converts the final interval length into an ME stationarity guarantee.

We use the class $\mathcal F_1(G,\rho,\Delta)$ and oracle model from Section~\ref{sec:introduction}. Set $Q:=\rho\Delta/\varepsilon^2$ and $\lambda:=1/(2\rho)$. Let $T_1^\star(G,\rho,\Delta,\varepsilon,P)$ denote the minimum worst-case number of rounds required to find a point satisfying $|f_\lambda'(x)|\le\varepsilon$ with probability at least $3/4$ on every instance, using at most $P$ queries per round.

\begin{theorem}[Optimal one-dimensional round complexity]
\label{thm:one-dimensional}
There exist universal constants $c_1>0$ and $C_1>1$ such that, if $0<\varepsilon\le c_1G$ and $Q\ge C_1$, then, for every integer $P\ge1$,
\[
T_1^\star(G,\rho,\Delta,\varepsilon,P)
=
\Theta\!\left(1+\frac{\log Q}{\log(P+1)}\right),
\]
with universal implicit constants. A deterministic algorithm using only exact function values attains the upper bound. The lower bound holds for randomized algorithms receiving exact function values and entire subdifferentials.
\end{theorem}

The same conclusion holds if $3/4$ in the definition of $T_1^\star$ is replaced by any $p\in(1/2,1]$, with implicit constants that can be chosen uniformly over all such $p$.

\subsection{Interval search using function values}
\label{subsec:one-dimensional-upper}

Set $\alpha:=\varepsilon/4$, $\delta:=\varepsilon/(4\rho)$, and $g(x):=f(x)+\alpha|x|$. The regularization serves two purposes: it converts the initial objective gap into a finite search interval, and it ensures that an interior local minimizer of $g$ has a small subgradient of $f$. The algorithm maintains three points $a<b<c$ satisfying $g(b)\le\min\{g(a),g(c)\}$. This condition guarantees that $g$ has an interior minimizer on $[a,c]$. We shrink the interval while preserving this condition and stop once its length is small enough to guarantee approximate stationarity.

To initialize, let $R:=\max\{\delta,\Delta/\alpha\}$. The gap bound gives $g(\pm R)\ge f(0)-\Delta+\alpha R\ge g(0)$. Query $f(0)$ and set $(a,b,c):=(-R,0,R)$. The endpoint values need not be queried, since the required inequalities follow from the gap bound.

For $P\ge2$, query $t_j:=a+j(c-a)/(P+1)$, $j=1,\ldots,P$, in parallel. Choose $b'$ to minimize $g$ among the old middle point $b$ and the new query points. Retaining $b$ among the candidates ensures $g(b')\le g(b)$, so the inequalities involving the old endpoints remain valid. Sort the distinct points in $\{a,b,t_1,\ldots,t_P,c\}$, and let $a'$ and $c'$ be the immediate neighbors of $b'$. Then $g(b')\le\min\{g(a'),g(c')\}$ and $c'-a'\le2(c-a)/(P+1)$. Replacing the current triple by $(a',b',c')$ therefore preserves the invariant and contracts the interval.

For $P=1$, query the midpoint of the longer side. If $c-b\ge b-a$, set $t:=(b+c)/2$ and replace $(a,b,c)$ by $(b,t,c)$ when $g(t)<g(b)$, and by $(a,b,t)$ otherwise. Treat the left side symmetrically. Initially, the two sides have equal length, and each update keeps the ratio of the longer side to the shorter side at most $2$. Indeed, writing $\ell:=b-a$ and $r:=c-b$, suppose $\ell\le r\le2\ell$. The two possible updates give side lengths $(r/2,r/2)$ and $(\ell,r/2)$, respectively. Both preserve the ratio bound, and their total lengths are at most $2(\ell+r)/3$ and $3(\ell+r)/4$. Thus, each round contracts the interval length by a factor of at most $3/4$.

Stop when $c-a\le\delta$ and return $b$. Continuity and the invariant imply that $g$ has a minimizer $z\in(a,c)$ over $[a,c]$. If the minimum is below $g(b)$, neither endpoint can attain it; otherwise, $b$ itself is an interior minimizer. Hence $0\in\partial g(z)=\partial f(z)+\alpha\partial|\cdot|(z)$, which gives $\operatorname{dist}(0,\partial f(z))\le\alpha$.

For a $\rho$-WC function and $\lambda=1/(2\rho)$, the ME satisfies $|f_\lambda'(z)|\le2\operatorname{dist}(0,\partial f(z))$, and its derivative is $2\rho$-Lipschitz; see Appendix~\ref{app:moreau-properties}. Consequently,
\[
|f_\lambda'(b)|
\le |f_\lambda'(z)|+2\rho|b-z|
\le2\alpha+2\rho\delta
=\varepsilon.
\]

The initial-to-final length ratio is $2R/\delta=\max\{2,32Q\}$. Define $\theta_1:=3/4$ and $\theta_P:=2/(P+1)$ for $P\ge2$. Including the initial query at the origin, the number of rounds satisfies
\[
T_{\mathrm{up}}
\le
1+\left\lceil
\frac{\log\max\{2,32Q\}}{\log(1/\theta_P)}
\right\rceil.
\]
For $P\ge2$, we have $\log((P+1)/2)\ge\frac13\log(P+1)$. Together with the contraction factor for $P=1$, this proves the upper bound in Theorem~\ref{thm:one-dimensional}. The algorithm uses only function values and does not require knowledge of $G$.

\subsection{Hard instances on nested intervals}
\label{subsec:one-dimensional-instances}

For the lower bound, we encode binary strings by nested intervals to hide the location of approximate stationary points. Such constructions also appear in the one-dimensional randomized lower bounds of \citet{KornowskiShamir2022}. Each bit selects one of two child intervals, and subsequent changes to the function are confined to the selected child.

A construction of depth $n$ will produce $2^n$ disjoint sets of approximate stationary points, with an initial objective gap of order $4^n\varepsilon^2/\rho$. The gap budget therefore permits $n=\Theta(\log Q)$, and locating an approximate stationary point requires information about the corresponding hidden interval.

Fix an integer $n\ge1$, and set $B:=40\varepsilon$ and $L_j:=16\cdot4^{n-j}B/\rho$ for $j=0,\ldots,n$. For each $\omega=(\omega_1,\ldots,\omega_n)\in\{0,1\}^n$, construct nested closed intervals $I_0\supset I_1\supset\cdots\supset I_n$, starting from $I_0:=[0,L_0]$. If $I_j=[a,a+L_j]$, define
\[
I_{j+1}:=
\begin{cases}
[a+L_j/8,\ a+3L_j/8], & \omega_{j+1}=0,\\
[a+5L_j/8,\ a+7L_j/8], & \omega_{j+1}=1.
\end{cases}
\]
The two candidate intervals are disjoint and have length $L_{j+1}=L_j/4$.

Assign the value $h_j$ to both endpoints of $I_j$, where $h_0:=0$ and $h_{j+1}:=h_j-BL_j/16$. For a parent interval $I_j=[a,c]$ and its selected child $I_{j+1}=[a',c']$, use continuously differentiable, piecewise quadratic functions on $[a,a']$ and $[c',c]$ to match the prescribed endpoint values. Their derivatives satisfy $-B\le f_\omega'(x)\le-B/20$ on $(a,a')$ and $B/20\le f_\omega'(x)\le B$ on $(c',c)$. Their smooth pieces have second derivative $0$ or $-\rho$, and the one-sided derivatives at the child endpoints satisfy $f_{\omega,-}'(a')=-B$ and $f_{\omega,+}'(c')=B$. Each connecting function depends only on the parent and selected child, not on later bits. Every nonterminal interval has length at least $64B/\rho$, so these conditions can be realized by linear and quadratic pieces. Explicit formulas are given in Appendix~\ref{app:interval-connectors}.

Write $I_n=[a_n,c_n]$ and define $m_\omega:=(a_n+c_n)/2$. On $I_n$, set $f_\omega(x):=h_n+B(|x-m_\omega|-L_n/2)$. Outside $I_0$, set $f_\omega(x):=-Bx$ for $x\le0$ and $f_\omega(x):=B(x-L_0)$ for $x\ge L_0$. This defines $f_\omega$ on the whole real line.

\begin{lemma}[Properties of the hard instances]
\label{lem:interval-instance-properties}
Each $f_\omega$ is $\rho$-WC and $B$-Lipschitz, with unique minimizer $m_\omega$ and objective gap $f_\omega(0)-\inf f_\omega=(4/3)(4^n+5)B^2/\rho$. For every $x\ne m_\omega$, we have $\operatorname{dist}(0,\partial f_\omega(x))\ge B/20=2\varepsilon$, while $\partial f_\omega(m_\omega)=[-B,B]$. Moreover, the set of approximate stationary points is exactly
\[
S_\omega
:=
\{x:|(f_\omega)_\lambda'(x)|\le\varepsilon\}
=
\left[m_\omega-\frac{\varepsilon}{2\rho},
      m_\omega+\frac{\varepsilon}{2\rho}\right].
\]
Each $S_\omega$ lies strictly inside its terminal interval $I_n$, and these sets are pairwise disjoint.
\end{lemma}

\begin{proof}
The connecting pieces have curvature at least $-\rho$, and all derivative jumps are upward, so $f_\omega+(\rho/2)|\cdot|^2$ is convex. The slopes have magnitude at most $B$. The function is strictly decreasing to the left of $m_\omega$ and strictly increasing to its right, with subgradient magnitude at least $B/20$ away from $m_\omega$. The verification at the joining points is included in Appendix~\ref{app:interval-connectors}.

The recursion for the endpoint values gives $h_n=-(4/3)(4^n-1)B^2/\rho$. Since $f_\omega(m_\omega)=h_n-BL_n/2$ and $BL_n/2=8B^2/\rho$, the stated objective gap follows.

To identify the stationary set, let $y:=\operatorname{prox}_{\lambda f_\omega}(x)$. Proximal optimality gives $(f_\omega)_\lambda'(x)=(x-y)/\lambda\in\partial f_\omega(y)$. If this derivative has magnitude at most $\varepsilon$, the subgradient bound forces $y=m_\omega$, and hence $|x-m_\omega|\le\lambda\varepsilon$. Conversely, this distance condition implies $(x-m_\omega)/\lambda\in[-\varepsilon,\varepsilon]\subseteq\partial f_\omega(m_\omega)$. Strong convexity of the proximal subproblem then gives $y=m_\omega$, proving the formula for $S_\omega$. Since $\lambda\varepsilon<L_n/2$ and distinct binary strings yield disjoint terminal intervals, the remaining claims follow.
\end{proof}

The construction also has an oracle locality property: once the first $j$ bits are fixed, both the function value and the entire subdifferential are determined at every point in $\mathbb R\setminus I_j$, independently of the remaining bits.

\subsection{Information bounds for randomized parallel queries}
\label{subsec:one-dimensional-lower}

Draw $\omega$ uniformly from $\{0,1\}^n$, independently of the algorithm's internal randomness. We analyze an augmented oracle that reveals a prefix of $\omega$ in addition to the original oracle replies.

Suppose the prefix known at the start of a round is $(\omega_1,\ldots,\omega_k)$, with corresponding interval $I_k$. All query points in the batch are chosen before any replies are received. For each query point $x$, start from this same prefix and reveal successive bits while $x$ belongs to the current closed interval and unrevealed bits remain. Each revealed bit selects the next child interval. The process stops when $x$ lies outside the current interval or all $n$ bits have been revealed.

For the whole batch, the oracle returns the longest of the prefixes required by the individual queries. This prefix contains all the others and determines the exact function values and entire subdifferentials at every query point. The algorithm can use it in subsequent rounds. Closed intervals account for endpoint queries, whose subdifferentials may depend on the next connecting function. Since the augmented oracle provides more information than the original oracle, a lower bound for the augmented model also applies to the original model.

\begin{lemma}[Information per round and success probability]
\label{lem:prefix-information}
Every randomized algorithm using the augmented oracle with at most $P$ queries per round and at most $T$ rounds on every execution satisfies
\[
\Pr_{\omega,\mathcal A}\{\widehat x\in S_\omega\}
\le
\min\left\{1,\frac{T(\log_2P+2)+2}{n}\right\}.
\]
A zero-round algorithm has average success probability at most $2^{-n}$.
\end{lemma}

\begin{proof}
Let $\mathcal H$ denote the history at the start of a round, including the algorithm's random seed and the revealed prefix. Conditional on $\mathcal H$, all query points in the current batch are fixed. Each revelation and stopping decision depends only on the query points and the bits already revealed, so the unrevealed bits remain independent and uniform.

For a fixed query point $x$, let $J_x$ be the number of additional bits required, measured from the prefix at the start of the round. To reveal at least $j$ bits, the first $j-1$ child selections must all choose an interval containing $x$. Since the two candidate children are disjoint, $\Pr(J_x\ge j\mid\mathcal H)\le2^{-(j-1)}$ for every $j\ge1$.

The additional prefix length for the batch is $J:=\max_x J_x$, where the maximum is over its query points. A union bound gives $\Pr(J\ge j\mid\mathcal H)\le\min\{1,P\,2^{-(j-1)}\}$. The tail-sum formula therefore yields
\[
\mathbb E[J\mid\mathcal H]
\le
\sum_{j\ge1}\min\{1,P\,2^{-(j-1)}\}
\le\log_2P+2.
\]

The output need not be a previously queried point. To account for this, perform one additional prefix query at $\widehat x$ after termination, solely for the analysis. This query is not counted as an algorithmic round. The single-query tail bound shows that it reveals at most $2$ additional bits in conditional expectation. If the algorithm succeeds, Lemma~\ref{lem:interval-instance-properties} places $\widehat x$ strictly inside the terminal interval, so all $n$ bits have been revealed after this additional query.

Let $Z$ denote the total number of revealed bits after the additional query. Averaging the conditional bounds over the histories and the algorithm's random seed gives
\[
n\,\Pr_{\omega,\mathcal A}\{\widehat x\in S_\omega\}
\le\mathbb E_{\omega,\mathcal A}Z
\le T(\log_2P+2)+2.
\]
Together with the trivial upper bound of $1$ on success probability, this proves the first claim.

With zero rounds, the output is independent of $\omega$. Since the $2^n$ sets $S_\omega$ are pairwise disjoint, any fixed output belongs to at most one of them. Averaging over the output gives success probability at most $2^{-n}$.
\end{proof}

\begin{proof}[Completion of the proof of Theorem~\ref{thm:one-dimensional}]
The upper bound was established in Section~\ref{subsec:one-dimensional-upper}. For the lower bound, we may take $c_1:=1/40$ and $C_1:=19200$. Choose $n:=\lfloor\log_4(3Q/6400-5)\rfloor$. Since $Q\ge19200$, we have $n\ge1$ and $n=\Theta(\log Q)$. With $B=40\varepsilon$, Lemma~\ref{lem:interval-instance-properties} gives
\[
f_\omega(0)-\inf f_\omega
=
\frac{6400}{3}(4^n+5)\frac{\varepsilon^2}{\rho}
\le\Delta.
\]
Also, $B\le G$, so every constructed function belongs to $\mathcal F_1(G,\rho,\Delta)$.

Fix any $p\in(1/2,1]$. If an algorithm succeeds with probability at least $p$ on every instance, its average success probability under the uniform distribution on $\omega$ is at least $p$. Lemma~\ref{lem:prefix-information} implies $T\ge(pn-2)/(\log_2P+2)$. For $n\ge8$, we have $pn-2\ge n/4$ and $\log_2P+2\le2\log_2(P+1)$, giving $T\ge n/(8\log_2(P+1))$.

The zero-round bound also gives $T\ge1$, since $2^{-n}\le1/2<p$. This constant lower bound covers $1\le n<8$ as well. Consequently,
\[
T
=
\Omega\!\left(1+\frac{n}{\log(P+1)}\right)
=
\Omega\!\left(1+\frac{\log Q}{\log(P+1)}\right).
\]
The constants are uniform over $p\in(1/2,1]$. The upper bound is deterministic and thus has success probability $1$, proving both the theorem and its extension to these success probabilities.
\end{proof}

These bounds show how the order structure of the real line allows function-value comparisons to shrink a search interval, with parallel queries accelerating this reduction. Regularization and weak convexity convert the interval length into an ME stationarity guarantee. The matching lower bound shows that this function-value search attains the optimal worst-case round complexity up to universal constants, even among randomized algorithms with access to the entire subdifferential.

\section{Conclusion}
\label{sec:conclusion}

We show that the standard subgradient rate remains optimal for randomized parallel algorithms in sufficiently high dimensions within the stated parameter regime. In one dimension, the order structure of the real line permits faster search, and we characterize precisely the reduction in round complexity provided by parallel queries.
A natural open problem is to characterize the dimension threshold at which the subgradient rate becomes unavoidable. For a prescribed accuracy, initial objective gap, and parallel query budget, this requires determining how large the dimension must be to sustain the matching lower bound, and what faster rates are attainable below that threshold. Such a characterization would clarify how low-dimensional geometry and parallel queries jointly affect the complexity of stationarity search.

\bibliographystyle{plainnat}
\bibliography{weak_convex_references}

@article{Vial1983,
  author  = {Vial, Jean-Philippe},
  title   = {Strong and Weak Convexity of Sets and Functions},
  journal = {Mathematics of Operations Research},
  year    = {1983},
  volume  = {8},
  number  = {2},
  pages   = {231--259},
  doi     = {10.1287/moor.8.2.231}
}

@article{DuchiRuan2018,
  author  = {Duchi, John C. and Ruan, Feng},
  title   = {Stochastic Methods for Composite and Weakly Convex Optimization Problems},
  journal = {SIAM Journal on Optimization},
  year    = {2018},
  volume  = {28},
  number  = {4},
  pages   = {3229--3259},
  doi     = {10.1137/17M1135086}
}

@article{DavisGrimmer2019,
  author  = {Davis, Damek and Grimmer, Benjamin},
  title   = {Proximally Guided Stochastic Subgradient Method for Nonsmooth, Nonconvex Problems},
  journal = {SIAM Journal on Optimization},
  year    = {2019},
  volume  = {29},
  number  = {3},
  pages   = {1908--1930},
  doi     = {10.1137/17M1151031}
}

@misc{DavisDrusvyatskiy2018,
  author        = {Davis, Damek and Drusvyatskiy, Dmitriy},
  title         = {Stochastic Subgradient Method Converges at the Rate {$O(k^{-1/4})$} on Weakly Convex Functions},
  year          = {2018},
  howpublished  = {arXiv preprint arXiv:1802.02988},
  eprint        = {1802.02988},
  archivePrefix = {arXiv},
  primaryClass  = {math.OC},
  url           = {https://arxiv.org/abs/1802.02988}
}

@article{DavisDrusvyatskiy2019,
  author  = {Davis, Damek and Drusvyatskiy, Dmitriy},
  title   = {Stochastic Model-Based Minimization of Weakly Convex Functions},
  journal = {SIAM Journal on Optimization},
  year    = {2019},
  volume  = {29},
  number  = {1},
  pages   = {207--239},
  doi     = {10.1137/18M1178244}
}

@misc{LiPan2026,
  author        = {Li, Jiajin and Pan, Siyu},
  title         = {Optimal Deterministic Oracle Complexity for Weakly Convex Optimization},
  year          = {2026},
  howpublished  = {arXiv preprint arXiv:2608.03246},
  eprint        = {2608.03246},
  archivePrefix = {arXiv},
  primaryClass  = {math.OC},
  url           = {https://arxiv.org/abs/2608.03246}
}

@article{Nemirovski1994,
  author  = {Nemirovski, Arkadi},
  title   = {On Parallel Complexity of Nonsmooth Convex Optimization},
  journal = {Journal of Complexity},
  year    = {1994},
  volume  = {10},
  number  = {4},
  pages   = {451--463},
  doi     = {10.1006/jcom.1994.1025}
}

@article{DiakonikolasGuzman2020,
  author  = {Diakonikolas, Jelena and Guzm{\'a}n, Crist{\'o}bal},
  title   = {Lower Bounds for Parallel and Randomized Convex Optimization},
  journal = {Journal of Machine Learning Research},
  year    = {2020},
  volume  = {21},
  number  = {5},
  pages   = {1--31},
  url     = {https://www.jmlr.org/papers/v21/19-771.html}
}

@article{CarmonEtAl2020,
  author  = {Carmon, Yair and Duchi, John C. and Hinder, Oliver and Sidford, Aaron},
  title   = {Lower Bounds for Finding Stationary Points {I}},
  journal = {Mathematical Programming},
  year    = {2020},
  volume  = {184},
  number  = {1--2},
  pages   = {71--120},
  doi     = {10.1007/s10107-019-01406-y}
}

@inproceedings{ZhouEtAl2025,
  author    = {Zhou, Huanjian and Han, Andi and Takeda, Akiko and Sugiyama, Masashi},
  title     = {The Adaptive Complexity of Finding a Stationary Point},
  booktitle = {Proceedings of the Thirty-Eighth Conference on Learning Theory},
  series    = {Proceedings of Machine Learning Research},
  year      = {2025},
  volume    = {291},
  pages     = {6091--6123},
  publisher = {PMLR},
  url       = {https://proceedings.mlr.press/v291/huanjian25a.html}
}

@inproceedings{BubeckMikulincer2020,
  author    = {Bubeck, S{\'e}bastien and Mikulincer, Dan},
  title     = {How to Trap a Gradient Flow},
  booktitle = {Proceedings of the Thirty-Third Conference on Learning Theory},
  series    = {Proceedings of Machine Learning Research},
  year      = {2020},
  volume    = {125},
  pages     = {940--960},
  publisher = {PMLR},
  url       = {https://proceedings.mlr.press/v125/bubeck20b.html}
}

@inproceedings{ChewiEtAl2023,
  author    = {Chewi, Sinho and Bubeck, S{\'e}bastien and Salim, Adil},
  title     = {On the Complexity of Finding Stationary Points of Smooth Functions in One Dimension},
  booktitle = {Proceedings of the 34th International Conference on Algorithmic Learning Theory},
  series    = {Proceedings of Machine Learning Research},
  year      = {2023},
  volume    = {201},
  pages     = {358--374},
  publisher = {PMLR},
  url       = {https://proceedings.mlr.press/v201/chewi23a.html}
}

@misc{ZhouSuzuki2026,
  author        = {Zhou, Huanjian and Suzuki, Taiji},
  title         = {On the Complexity of Finding Stationary Points of {$p$}-th Order {H{\"o}lder} Smooth Functions in One Dimension},
  year          = {2026},
  howpublished  = {arXiv preprint arXiv:2610.16295},
  eprint        = {2610.16295},
  archivePrefix = {arXiv},
  primaryClass  = {math.OC},
  url           = {https://arxiv.org/abs/2610.16295}
}

@inproceedings{DavisEtAl2022,
  author    = {Davis, Damek and Drusvyatskiy, Dmitriy and
               Lee, Yin Tat and Padmanabhan, Swati and Ye, Guanghao},
  title     = {A Gradient Sampling Method with Complexity Guarantees
               for {Lipschitz} Functions in High and Low Dimensions},
  booktitle = {Advances in Neural Information Processing Systems},
  volume    = {35},
  pages     = {6692--6703},
  publisher = {Curran Associates, Inc.},
  year      = {2022},
  doi       = {10.52202/068431-0485},
  url       = {https://proceedings.neurips.cc/paper_files/paper/2022/hash/2c8d9636f74d0207ff4f65956010f450-Abstract-Conference.html}
}

@misc{KornowskiShamir2022,
  author        = {Kornowski, Guy and Shamir, Ohad},
  title         = {On the Complexity of Finding Small Subgradients
                   in Nonsmooth Optimization},
  year          = {2022},
  eprint        = {2209.10346},
  archivePrefix = {arXiv},
  primaryClass  = {math.OC},
  note          = {Revised April 2025},
  url           = {https://arxiv.org/abs/2209.10346v2}
}

\appendix

\section{Moreau-envelope estimates}
\label{app:moreau-properties}

Let $f:\mathbb R^d\to\mathbb R$ be finite-valued and $\rho$-WC, and set $\lambda:=1/(2\rho)$. Write $p(x):=\operatorname{prox}_{\lambda f}(x)$ and $v(x):=\nabla f_\lambda(x)$. The proximal objective $f(y)+\rho\|y-x\|^2$ is strongly convex and coercive, so $p(x)$ is uniquely defined. We establish the two estimates used in the one-dimensional upper bound:
\[
\|v(x)\|\le2\operatorname{dist}(0,\partial f(x)),
~~
\|v(x)-v(y)\|\le2\rho\|x-y\|,
\]
where $\operatorname{dist}(0,C):=\inf_{g\in C}\|g\|$.

\begin{proof}[Proof of the envelope estimates]
Convexity of $f+(\rho/2)\|\cdot\|^2$ implies
\[
\langle g-h,x-y\rangle\ge-\rho\|x-y\|^2
~\text{for }g\in\partial f(x),\ h\in\partial f(y).
\]
Proximal optimality gives $2\rho(x-p(x))\in\partial f(p(x))$. Applying the preceding inequality at $p(x)$ and $p(y)$ yields $\langle x-y,p(x)-p(y)\rangle\ge\frac12\|p(x)-p(y)\|^2$. In particular, $p$ is $2$-Lipschitz. The unique-minimizer differentiation formula therefore gives $v(x)=2\rho(x-p(x))$, and
\[
\|v(x)-v(y)\|^2
=4\rho^2\|(x-y)-(p(x)-p(y))\|^2
\le4\rho^2\|x-y\|^2.
\]
Next, for any $g\in\partial f(x)$, apply the same subgradient inequality at $x$ and $p(x)$. Substituting $x-p(x)=v(x)/(2\rho)$ gives $\frac12\|v(x)\|^2\le\langle g,v(x)\rangle\le\|g\|\|v(x)\|$. Taking the infimum over $g$ proves the first estimate, including the case $v(x)=0$.
\end{proof}

The same identities give the subgradient upper bound used in Corollary~\ref{cor:high-dimensional-optimal}. Suppose also that $f$ is $G$-Lipschitz and $f(x_0)-\inf f\le\Delta$. For $g_t\in\partial f(x_t)$ and $x_{t+1}:=x_t-\eta g_t$, let $p_t:=p(x_t)$. Comparing the proximal objective at $p_t$ and $x_t$, and then using weak convexity at $x_t$, gives $\langle g_t,x_t-p_t\rangle\ge(\rho/2)\|x_t-p_t\|^2$. Evaluating the next envelope at $p_t$ consequently yields
\[
\begin{aligned}
f_\lambda(x_{t+1})
&\le f_\lambda(x_t)-2\rho\eta\langle g_t,x_t-p_t\rangle+\rho\eta^2\|g_t\|^2\\
&\le f_\lambda(x_t)-\frac\eta4\|v(x_t)\|^2+\rho\eta^2G^2.
\end{aligned}
\]
Since $f_\lambda(x_0)-\inf f\le\Delta$, summing over $t=0,\ldots,T-1$ and choosing $R$ uniformly from these indices gives $\mathbb E\|v(x_R)\|^2\le4\Delta/(\eta T)+4\rho\eta G^2$. The choice $\eta:=\sqrt{\Delta/(\rho G^2T)}$ gives the bound $8\sqrt{\rho G^2\Delta/T}$ stated in the main text.

\section{Properties of the base function}
\label{app:base-properties}

\begin{proof}[Proof of Proposition~\ref{prop:base-properties}]
Write $h:=h_{M,N}$ and $\phi_r:=\phi(u^{(r)})$. Let $\Lambda:=\{\lambda\in\mathbb R^M:\lambda\ge0,\ \|\lambda\|\le1\}$ and set $\lambda_{M+1}:=0$. The support-function identity $\|[w]_+\|=\max_{\lambda\in\Lambda}\langle\lambda,w\rangle$ gives $h(u)=\max_{\lambda\in\Lambda}h_\lambda(u)$, where
\[
h_\lambda(u):=\frac2N\left[
\sum_{r=1}^M\lambda_r\phi_r
-\sum_{r=1}^M(1-\lambda_{r+1})p_r
+\lambda_1-\frac54\sum_{r=1}^M\lambda_r
\right].
\]
The derivative $q'$ is globally $6$-Lipschitz and satisfies $0\le q'\le3/2$. Hence $\nabla\psi$ is $3N/8$-Lipschitz, with nonnegative coordinates bounded by $3/8$. Since $0\le1-\lambda_{r+1}\le1$, each $h_\lambda$ is $3/4$-WC: the affine maxima are convex, and the smooth terms act on disjoint blocks. Taking the maximum preserves this weak convexity parameter.

Let $\Lambda(u):=\arg\max_{\lambda\in\Lambda}\langle\lambda,\sigma(u)\rangle$. Applying the subdifferential rule for a compact maximum after adding the common convexifying quadratic shows that $\partial h(u)$ is the convex hull of the vectors $g$ with block components
\[
g^{(r)}=\frac2N\left[\lambda_r\xi_r-(1-\lambda_{r+1})\nabla\psi(u^{(r)})\right],
~~ r\in[M],
\]
where $\lambda\in\Lambda(u)$ and $\xi_r\in\partial\phi(u^{(r)})$. Every $\xi_r$ is a convex combination of the active vectors among $0,-(N/8)e_1,\ldots,-(N/8)e_N$, where $e_j$ denotes the $j$th coordinate vector. Thus all coordinates of every subgradient of $h$ are nonpositive. Also, $\|\xi_r\|\le N/8$ and $\|\nabla\psi(u^{(r)})\|\le3\sqrt N/8$, giving
\[
\|g\|\le\frac2N\left(\frac N8\|\lambda\|+\frac38\sqrt{MN}\right)
\le\frac14+\frac34\sqrt{\frac MN}<\frac12.
\]
The same bound holds for convex combinations and, equivalently, uniformly bounds the Lipschitz constants of the functions $h_\lambda$. It therefore bounds $\operatorname{Lip}(h)$.

If $S(u)>0$, the maximizing weight is unique: $\lambda=[\sigma(u)]_+/S(u)$, with $\|\lambda\|=1$. Whenever $\lambda_r>0$, we have $\sigma_r>0$ and $\phi_r>5/4-p_{r-1}\ge1/4$. The zero branch of $\phi_r$ is then inactive, so $\xi_r=-(N/8)w_r$ for some $w_r\ge0$ with $\sum_j(w_r)_j=1$. Consequently, $\|\xi_r\|\ge\sqrt N/8$. The two terms in each block of $g$ have the same coordinatewise sign, and hence
\[
\|g\|^2\ge\frac4{N^2}\sum_{r=1}^M\lambda_r^2\|\xi_r\|^2
\ge\frac1{16N}\sum_{r=1}^M\lambda_r^2
=\frac1{16N}.
\]
For this fixed $\lambda$, the displayed subgradient representation already defines a convex set, so the estimate holds for every $g\in\partial h(u)$.

Since $S\ge0$ and $0\le p_r\le1$, we have $h\ge-2M/N$. At the origin, $p_r=0$ for $r\ge1$, $\sigma_1=c_1-1/4>0$, and $\sigma_r=c_1-5/4<0$ for $r\ge2$. Thus $h(0)=2(c_1-1/4)/N\in[0,2/N)$, and $h(0)-\inf h\le(2M+2)/N\le3M/N$.

Finally, suppose $S(u)=0$. Since $p_0=1$, the inequality $\sigma_1\le0$ implies $\phi_1\le1/4$, and therefore $u^{(1)}_j\ge(8/N)(c_j-1/4)\ge6/N$ for every $j$. The arguments of $q$ are then at least $3/2-1/(16N)>1$, so $p_1=1$. Induction gives $p_r=1$ for every $r$. In particular, $p_M<27/32$ implies $S(u)>0$.
\end{proof}

\section{Exact locality of the base oracle}
\label{app:base-locality}

\begin{proof}[Proof of Lemma~\ref{lem:base-locality}]
Set $\bar a:=3a/2$. Let $\mathcal T_t$ be the open set of all $v$ satisfying $|v_j^{(r)}|<\bar a$ whenever $(r,j)\notin\mathcal J_{t-1}$. Both $u$ and $Q_tu$ belong to $\mathcal T_t$. We show that $h_{M,N}(v)=h_{M,N}(Q_tv)$ throughout this set.

For each block, let $k_r:=|\{j:(r,j)\in\mathcal J_{t-1}\}|$. If $j>k_r$, then $v_j^{(r)}<\bar a<2a$, so its contribution to the progress function and its derivative is zero. In particular, all progress values are preserved by $Q_t$.

Suppose $0<k_r<N$. The coordinate $k_r+1$ is retained by $Q_t$, and for every $j>k_r+1$,
\[
\left(c_{k_r+1}-\frac N8v_{k_r+1}^{(r)}\right)
-\left(c_j-\frac N8v_j^{(r)}\right)
\ge\frac1{16N}-\frac N4\bar a
=\frac1{64N}>0.
\]
Thus every deleted affine branch is strictly below a retained branch. The same inequality holds after projection, so $\phi(v^{(r)})$ is unchanged. A block with $k_r=N$ is retained in full. If $r=1$ and $k_1=0$, then $t=0$; the first coordinate is retained, and the same comparison shows that its affine branch strictly dominates all deleted branches.

It remains to consider $k_r=0$ with $r\ge2$. At most $L$ coordinates of block $r-1$ belong to $\mathcal J_{t-1}$, giving $p_{r-1}(v)\le L/N\le3/16$. All coordinates of block $r$ have magnitude below $\bar a$, and hence
\[
\sigma_r(v)\le c_1+\frac N8\bar a+\frac3{16}-\frac54
=-\frac5{128N}<0.
\]
This estimate also holds at $Q_tv$. Such a block contributes zero to both $[\sigma]_+$ and its own progress value. These cases prove the required identity on $\mathcal T_t$.

Local equality and the subdifferential rule for linear composition yield $\partial h_{M,N}(u)=Q_t\partial h_{M,N}(Q_tu)$. Applying the same identity at $Q_tu$ gives $\partial h_{M,N}(Q_tu)=Q_t\partial h_{M,N}(Q_tu)$. Combining the two equalities proves both the equality of the full subdifferentials and their inclusion in $\operatorname{range}Q_t$.
\end{proof}

\section{Properties of a single stage}
\label{app:single-stage}

\begin{proof}[Proof of Lemma~\ref{lem:single-stage}]
Write $R_N(u):=\frac14\|[-u-2a\mathbf1]_+\|$ and $s:=\langle e,u\rangle$, so that $b_N(u)=h_N(u)+R_N(u)-\tau_N(s)$. The function $R_N$ is convex and $1/4$-Lipschitz, whereas $0\le\tau_N'\le1/4$ and $\operatorname{Lip}(\tau_N')=1/20$. Proposition~\ref{prop:base-properties} therefore gives weak convexity parameter $3/4+1/20=4/5$ and Lipschitz constant at most $3/4+\sqrt3/8<31/32$. The subdifferential sum rule gives $g=g_h+g_R-\tau_N'(s)e$, with $g_h\in\partial h_N(u)$ and $g_R\in\partial R_N(u)$. All three summands are componentwise nonpositive, so the norm of their sum is at least the norm of each summand.

If $S(u)>0$, Proposition~\ref{prop:base-properties} gives $\|g_h\|\ge1/(4\sqrt N)$. If $S(u)=0$, the proof of Proposition~\ref{prop:base-properties} gives $u_j^{(r)}\ge6/N$ in every block. Thus $s\ge6/\sqrt N$ and $\tau_N'(s)\ge1/(4\sqrt N)$, proving the global subgradient bound in this case as well.

If $R_N(u)=0$, every coordinate is at least $-2a$, so $\phi(u^{(r)})\le17/16$ and $\sigma_r\le13/16<1$. Consequently, $b_N(u)\le h_N(u)\le2\sqrt m/N<1/4$. Hence $b_N(u)\ge1/4$ forces $R_N(u)>0$, where $R_N$ is differentiable and $\|\nabla R_N(u)\|=1/4$. Conversely, $h_N\ge-2m/N\ge-1/6$ and $R_N\ge0$ imply that $b_N(u)\le-1$ forces $\tau_N(s)\ge5/6$. Since $\tau_N(1/\sqrt N+5)=5/8$, this lies in the linear part of $\tau_N$, where $\tau_N'(s)=1/4$. The common coordinatewise sign proves $\|g\|\ge1/4$ outside $(-1,1/4)$.

It remains to establish the inner-product estimate. Use the subgradient representation from Appendix~\ref{app:base-properties}, with $\lambda\ge0$, $\|\lambda\|\le1$, $\lambda_{m+1}:=0$, and $\xi_r\in\partial\phi(u^{(r)})$. Each such representing vector $w$ has blocks $w^{(r)}=(2/N)[\lambda_r\xi_r-(1-\lambda_{r+1})\nabla\psi(u^{(r)})]$, and $\sum_r\lambda_r\sigma_r=S(u)$. The affine branches defining $\phi$ give $\langle\xi_r,u^{(r)}\rangle\ge\phi(u^{(r)})-17/16$. To bound the progress contribution, put $t_j:=\frac N4(u_j^{(r)}-2a)$. Since $q'$ vanishes outside $(0,1)$, $\max_{0\le t\le1}tq'(t)=8/9$, and $\max q'=3/2$, we have
\[
0\le\langle\nabla\psi(u^{(r)}),u^{(r)}\rangle
=\sum_{j=1}^N q'(t_j)\left(\frac{t_j}{N}+\frac a2\right)
\le\frac89+\frac{3}{32N}<1.
\]
Using $0\le p_r\le1$ and $\sum_r\lambda_r\le\sqrt m$, it follows that
\[
\begin{aligned}
\langle w,u\rangle-h_N(u)
\ge\frac2N\left[\sum_{r=1}^m\lambda_r\left(\frac3{16}-p_{r-1}\right)-m\right]
\ge-\frac{13\sqrt m}{8N}-\frac{2m}{N}.
\end{aligned}
\]
The same inequality holds for every $g_h\in\partial h_N(u)$ by convex combinations. Convexity and $R_N(0)=0$ give $\langle g_R,u\rangle\ge R_N(u)$. Finally, direct integration gives $\tau_N(s)-s\tau_N'(s)\ge-5/8-1/(4\sqrt N)$: writing $t:=(s-1/\sqrt N)_+$, the left side is $-t^2/40-t/(20\sqrt N)$ for $0\le t\le5$ and is constant for $t\ge5$. Therefore,
\[
\langle g,u\rangle-b_N(u)
\ge-\frac{13\sqrt m}{8N}-\frac{2m}{N}-\frac58-\frac1{4\sqrt N}
>-1>-2,
\]
where $N\ge120$ and $m\le N/12$ were used.

At the origin, both added terms vanish, giving $b_N(0)=h_N(0)\in[0,2/N)$. If $s\le1/\sqrt N$, the terminal term vanishes and $b_N(u)\ge h_N(u)\ge-1/6$. Along $u+te$ as $t\to\infty$, the last-block affine and progress terms and the penalty eventually become constant, while $\tau_N(\langle e,u\rangle+t)$ grows linearly with slope $1/4$. Hence $b_N(u+te)\to-\infty$.
\end{proof}

\section{Radial extension of the combined function}
\label{app:radial-properties}

\begin{proof}[Remaining statements of Proposition~\ref{prop:radial-extension}]
Composition with $U^\top$ preserves the Lipschitz and weak convexity bounds of $F_K$, since $\|U^\top\|=1$. The radial branch is convex and $1$-Lipschitz, so their maximum is $4/5$-WC and $1$-Lipschitz. At the origin, $-A<0\le F_K(0)=h_N(0)$, giving $\Phi_U(0)=h_N(0)$.

To show that the lower bound $-2K$ is attained, set $\bar z^{(i)}:=(17/(2N))\mathbf 1_{d_0}+(8K+5)e$ for $i\in[K]$. Every coordinate is at least $17/(2N)$, so all affine maxima vanish, all progress values equal $1$, and the negative-coordinate penalties vanish. Thus $h_N(\bar z^{(i)})=-2m/N$. Since $\langle e,\bar z^{(i)}\rangle-1/\sqrt N=8K+5+15/(2\sqrt N)>5$, the linear part of $\tau_N$ gives
\[
\tau_N(\langle e,\bar z^{(i)}\rangle)
=\frac14\left(8K+5+\frac{15}{2\sqrt N}\right)-\frac58
=2K+\frac58+\frac{15}{8\sqrt N}>2K.
\]
Consequently, $F_K(\bar z)=-2K$. Moreover,
\[
\|\bar z\|
\le\sqrt K\left(\frac{17}{2}\sqrt{\frac mN}+8K+5\right)
<16K\sqrt K
<16(2K+2)\sqrt{N+K-1}=A-2K.
\]
At $x:=U\bar z$, the radial branch is therefore below $-2K$. It follows that $\inf\Phi_U=-2K$ and $\Phi_U(0)-\inf\Phi_U=2K+h_N(0)\in[2K,3K]$.

For any coordinate projection $\Pi$, including the identity, the Lipschitz bound gives $F_K(\Pi U^\top x)\le h_N(0)+(31/32)\|x\|$. If $\|x\|\ge64(A+1)$, then
\[
\|x\|-A-F_K(\Pi U^\top x)
\ge\frac{\|x\|}{32}-A-h_N(0)
\ge A+2-h_N(0)>0.
\]
The strict inequality persists locally, so both the value and the entire subdifferential of the maximum agree with those of the radial branch. This proves the oracle identity for the full and truncated functions.

For completeness, the convex-combination estimate used in the main text follows directly from $\|e_x\|=1$ and $\langle v,e_x\rangle\ge-\|v\|/4$. For $0\le\lambda\le1$, put $a_0:=\lambda\|v\|$ and $b_0:=1-\lambda$. Then
\[
\|\lambda v+(1-\lambda)e_x\|^2
\ge a_0^2+b_0^2-\frac12a_0b_0
\ge\frac34(a_0^2+b_0^2)
\ge\frac38(a_0+b_0)^2
\ge\frac38\min\{\|v\|,1\}^2.
\]
\end{proof}

\section{Auxiliary estimates for the small-gap construction}
\label{app:small-gap-construction}

\begin{proof}[Auxiliary estimates for Proposition~\ref{prop:normalized-small-gap}]
Write $h:=h_{M,N}$ and retain $s:=\sqrt{M/N}$, $A_s:=34s+2s^2$, and $B_s:=73s$. We first bound the radial component of a base-function subgradient. If $S(z)>0$ and $g\in\partial h(z)$, the subdifferential formula in Appendix~\ref{app:base-properties} expresses $g$ as a sum of componentwise nonpositive contributions from active affine branches and progress derivatives. An active affine branch with a nonzero coefficient satisfies $z_j^{(r)}\le8c_j/N<17/(2N)$. A nonzero progress derivative requires $2a<z_j^{(r)}<2a+4/N<17/(2N)$. Hence $g_{r,j}<0$ implies $z_j^{(r)}<17/(2N)$, and
\[
\langle g,z\rangle
\ge-\frac{17}{2N}\|g\|_1
\ge-\frac{17}{2}\sqrt{\frac MN}\,\|g\|
=-\frac{17}{2}s\|g\|.
\]

If the radial branch is active at $x$, the lower bound $h\ge-2s^2$ gives $\|x\|\ge A_s-2s^2=34s$. When both branches are active and $S(U^\top x)>0$, every $g\in\partial h(U^\top x)$ therefore satisfies $\langle Ug,x/\|x\|\rangle\ge-\|g\|/4$. The estimate in Appendix~\ref{app:radial-properties} and Proposition~\ref{prop:base-properties} imply
\[
\|\lambda Ug+(1-\lambda)x/\|x\|\|
\ge\sqrt{\frac38}\min\{\|g\|,1\}
\ge\frac{\sqrt{3/8}}{4\sqrt N}
>\frac1{8\sqrt N}
~(0\le\lambda\le1).
\]
The same lower bound holds when only the base branch or only the radial branch is active. Thus $S(U^\top x)>0$ implies $\operatorname{dist}(0,\partial\Psi_U(x))>1/(8\sqrt N)$.

Next, $s\le1/\sqrt{12}<1/2$ gives $A_s\le35s$, while $h(0)<2/N<s$ and $\operatorname{Lip}(h)<1/2$. For any coordinate projection $Q$, including the identity, and any $\|x\|\ge72s$, we consequently have
\[
h(QU^\top x)
<s+\frac12\|x\|
\le\|x\|-35s
\le\|x\|-A_s.
\]
The radial branch strictly dominates on a neighborhood of every such point. The full and truncated oracle replies are therefore $(\|x\|-A_s,\{x/\|x\|\})$.

Finally, suppose $|\langle u_{M,j},x\rangle|\le a$ for every $j>\lfloor N/2\rfloor$, and assume that $v:=\nabla(\Psi_U)_{1/2}(x)$ satisfies $\|v\|\le1/(16\sqrt N)$. Let $y:=\operatorname{prox}_{\Psi_U/2}(x)=x-v/2$, so that $v\in\partial\Psi_U(y)$. For a constrained last-block coordinate, $\langle u_{M,j},y\rangle-2a\le|\langle u_{M,j},v\rangle|/2$. Using $q(t)\le3(t_+)^2$ and the orthonormality of the columns of $U$ yields
\[
\begin{aligned}
\psi((U^\top y)^{(M)})
\le\frac{\lfloor N/2\rfloor}{N}
+\frac{3N}{64}\sum_{j>\lfloor N/2\rfloor}|\langle u_{M,j},v\rangle|^2
\le\frac12+\frac{3N}{64}\|v\|^2
\le\frac12+\frac3{16384}<\frac{27}{32}.
\end{aligned}
\]
Proposition~\ref{prop:base-properties} gives $S(U^\top y)>0$, so the subgradient bound above contradicts $v\in\partial\Psi_U(y)$. This proves the claimed ME gradient lower bound under the last-block condition.

If instead $\|x\|\ge B_s$, the same small-gradient assumption gives $\|x-y\|\le1/(32\sqrt N)<s$, and hence $\|y\|>72s$. The radial oracle identity then gives $\|v\|=1$, again a contradiction. These estimates establish both geometric properties used in the coupling argument of Proposition~\ref{prop:normalized-small-gap}.
\end{proof}

\section{Explicit connectors for the one-dimensional construction}
\label{app:interval-connectors}

We give the connecting functions used in Section~\ref{subsec:one-dimensional-instances} and verify their behavior at the joining points. Fix a parent interval $[a,c]$ of length $L\ge64B/\rho$, with endpoint value $h$, and either of its two candidate children $[a',c']$. The left and right connector lengths $\ell_-:=a'-a$ and $\ell_+:=c-c'$ belong to $\{L/8,5L/8\}$. The prescribed drop in endpoint value is $\eta:=BL/16$.

For either connector length $\ell$, define
\[
\begin{aligned}
\mu_\ell&:=B-\rho\ell+
\sqrt{(\rho\ell-B)^2+2\rho\eta-B^2},
~~
w_\ell:=\frac{B-\mu_\ell}{\rho},\\
C_\ell(t)&:=-\mu_\ell t-
\frac{\rho}{2}(t-\ell+w_\ell)_+^2
~ (0\le t\le\ell).
\end{aligned}
\]
The left connector is $f_\omega(x):=h+C_{\ell_-}(x-a)$ for $a\le x\le a'$, and the right connector is $f_\omega(x):=h+C_{\ell_+}(c-x)$ for $c'\le x\le c$. These formulas depend only on the parent interval, its endpoint value, and the selected child.

\begin{proof}[Verification of the connector properties]
Consider $A_\ell(\mu):=\mu\ell+(B-\mu)^2/(2\rho)$ on $[B/20,B]$. Since $\ell\ge L/8\ge8B/\rho$, its derivative satisfies $A_\ell'(\mu)=\ell-(B-\mu)/\rho>0$. Moreover,
\[
A_\ell(B/20)
\le\frac{BL}{32}+\frac{B^2}{2\rho}
\le\frac{5BL}{128}
<\eta
<B\ell=A_\ell(B).
\]
Thus, $A_\ell(\mu)=\eta$ has a unique solution in $(B/20,B)$, and solving the quadratic gives the stated $\mu_\ell$. In particular, $0<w_\ell<B/\rho<\ell$. The function $C_\ell$ is continuously differentiable and piecewise quadratic, with $C_\ell(0)=0$, $C_\ell(\ell)=-\eta$, and $C_\ell'(t)=-\mu_\ell-\rho(t-\ell+w_\ell)_+$. Hence $-B\le C_\ell'(t)\le-\mu_\ell<-B/20$, its smooth pieces have second derivative $0$ or $-\rho$, and $C_\ell'(\ell)=-B$. Reflection gives the positive derivative bounds for the right connector and derivative $B$ at $c'$.

It remains to check the joins between different levels. At the left endpoint of any nonterminal interval, the incoming derivative from its parent connector is $-B$, whereas its own left connector starts with derivative $-\mu_{\ell_-}\ge-B$. At the right endpoint, its own right connector ends with derivative $\mu_{\ell_+}\le B$, whereas the outgoing derivative of the parent connector is $B$. Both jumps are therefore upward. The same statements hold at the endpoints of the root interval, where the exterior linear pieces have derivatives $-B$ and $B$. At the terminal interval endpoints, the derivatives of the V-shaped terminal function agree with those of the parent connectors. Its midpoint has the upward jump from $-B$ to $B$. The joins within each connector have continuous derivatives.

Consequently, $f_\omega$ is continuous, all its smooth pieces have curvature at least $-\rho$, and all derivative jumps are upward. The derivative of $f_\omega(x)+(\rho/2)x^2$ is therefore nondecreasing, including across the joins, which proves weak convexity. The derivative bounds also prove that $f_\omega$ is $B$-Lipschitz. At every point, its full subdifferential is the interval $\partial f_\omega(x)=[f'_{\omega,-}(x),f'_{\omega,+}(x)]$. At a nonterminal left endpoint this interval is $[-B,-\mu_{\ell_-}]$, and at a nonterminal right endpoint it is $[\mu_{\ell_+},B]$. At the terminal endpoints it is respectively $\{-B\}$ and $\{B\}$, while at the terminal midpoint it is $[-B,B]$. Thus, every subgradient away from the midpoint has magnitude at least $B/20$, including at all joining points.

Finally, fix a prefix of length $j$. All continuations have the same intervals, endpoint values, and connecting functions outside $I_j$. Every point in $\mathbb R\setminus I_j$ has a neighborhood disjoint from $I_j$, so all continuations agree on that neighborhood and have the same function value and full subdifferential there. At the endpoints of $I_j$, the inner one-sided derivatives may depend on the next child through $\mu_{\ell_-}$ or $\mu_{\ell_+}$. This explains the use of closed intervals in the prefix oracle: an endpoint query reveals the next bit, after which it lies outside the selected child and its complete reply is determined.
\end{proof}

\end{document}